\documentclass[a4paper]{amsart}
\date{February 28, 2016}

\def\rd{\mathrm{d}}
\def\gM{\mathfrak{M}}
\def\rz{\mathbb{R}}
\def\bS{\mathbb{S}}
\newtheorem{thm}{Theorem}

\theoremstyle{definition}

\theoremstyle{remark}
\newtheorem{rem}{Remark}

\begin{document}
\title[Keller-Segel Equation]{Blow-Up of Solutions to the
  Patlak-Keller-Segel Equation in Dimension $\nu\geq2$}
\author[L. Chen]{Li Chen} \address{Mathematisches
  Institut\\Universit\"at Mannheim\\A5 6\\68131 Mannheim\\Germany}
\email{chen@math.uni-mannheim.de} \author[H. Siedentop]{Heinz
  Siedentop} \address{Mathematisches
  Institut\\Ludwig-Maximilians-Universit\"at
  M\"unchen\\Theresienstra\ss e 39, 80333 M\"unchen\\Germany}
\email{h.s@lmu.de}
%\maketitle

\begin{abstract}
  We prove a blow-up criterion for the solutions to the
  $\nu$-dimensional Patlak-Keller-Segel equation in the whole
  space. The condition is new in dimension three and higher. In
  dimension two it is exactly Dolbeault's and Perthame's blow-up
  condition, i.e., blow-up occurs, if the total mass exceeds $8\pi$.
\end{abstract}
\maketitle
\section{Introduction}

It is well-known that Patlak-Keller-Segel type equations describing
chemotaxis (Patlak \cite{Patlak1953} and Keller and Segel
\cite{KellerSegel1970,KellerSegel1971}) allow for both diffusion and
aggregation phenomena: depending on the initial data, the solution
might exist globally in time or blow up in finite time. Since
J\"ager's and Luckhaus' \cite{JagerLuckhaus1992} pioneering work the
analysis of this system proliferated. It is too ambitious to mention
all the important results in this short article, instead we refer to
the review of Wang et al \cite{Wangetal2016}.

We will focus on the blow-up of solutions of the
following system
\begin{eqnarray}
\nonumber \partial_t\rho -\nabla\cdot (\nabla\rho-\rho \nabla c)=0, &&\\
\label{equ}-\Delta c=\rho-Z\delta, &&x\in\rz^\nu, t\geq 0,\\
\nonumber\rho(x,0)=\rho_0(x)\geq 0,&&
\end{eqnarray}
where $\delta$ is the Dirac delta function at the origin. The
parameter $Z\in\rz$ represents the strength of the external
source. The source is repulsive, if $Z$ is positive, and attractive,
if $Z$ is negative.

For sake of simplicity, we write
\begin{eqnarray*}
 \gM_\mu(t):=\int_{\rz^\nu}|x|^\mu\rho(x,t)\rd x, \qquad \mu=0,1,\cdots
\end{eqnarray*}
for the moments of a solution $\rho$ at time $t$.  The system
\eqref{equ} conserves the total mass, i.e., for all times $t$ for
which the solution exists
$$ \gM_0(t)= \gM_0(0)=:\gM_0.
$$

For $\nu=2$ and $Z=0$, Dolbeault and Perthame
\cite{DolbeaultPerthame2004} and -- later in greater detail --
Blanchet et al \cite{Blanchetetal2006} showed that, if the mass
$\gM_0$ exceeds $8\pi$, then any classical solution blows up in finite
time, while, if $\gM_0< 8\pi$, then a classical solution exists
globally, since the diffusion dominates the aggregation which follows
from the logarithmic Hardy-Littlewood-Sobolev inequality.

Also for $\nu=2$, Wolansky and Espejo \cite{WolanskyEspejo2013}, show
that adding a repulsive point source will slow down the blow-up compared to
$Z=0$, where the blow-up occurs for masses that exceed $8\pi$, while
adding an attractive point source will enhance blow-up. This
generalizes \cite{DolbeaultPerthame2004}.

In higher dimensions, i.e., $\nu>2$, it is known (Perthame
\cite[Chapter 6]{Perthame2007}, Chapter 6) that there exists a
constant $c$ such that, if the initial datum fulfills
$\|\rho_0\|_{L^{\nu/2}}\leq c$, then a solution exists globally. On
the other hand there exists a (small) constant $c$ such that, if the
initial datum fulfills
\begin{equation}
\label{m2m0}
\gM_2<c \gM_0^\frac{\nu}{\nu-2},
\end{equation}
then all classical solutions blow up in finite time. However, the
results are much weaker than the corresponding version in dimension
two. This paper is a contribution to fill this gap.

A tool often used to show blow-up is the time evolution of the second
moment $\gM_2(t)$: if one can show that the time derivative is
strictly negative, then the finite time blow-up is proved. However,
the second moment is -- in some sense -- only natural in dimension
two.

The blow-up of solutions has also been studied for other variants of
the Patlak-Keller-Segel system. An example is
$\partial_t\rho -\nabla\cdot (\nabla\rho^m-\rho \nabla c)=0$, where a
porous media type degenerate diffusion is included. Sugiyama
\cite{Sugiyama2006} and Blanchet et al \cite{Blanchetetal2006}
obtained both existence and blow-up results for $m=2-\frac{2}{\nu}$,
and Chen et al \cite{Chenetal2012} for
$m=\frac{2\nu}{\nu+2}$. Furthermore, Chen and Wang \cite{ChenWang2014}
obtained sharp results on the occurrence of the blow-up and existence
for $m$ in between $2-\frac{2}{\nu}$ and $\frac{2\nu}{\nu+2}$.  In
particular they were able to identify very large classes of initial
conditions for which they can predict blow-up in finite time
respectively existence of solutions. Another variant is the model
with nonlinear chemotaxis sensitivity
$\partial_t\rho -\nabla\cdot (\nabla\rho-\rho^\alpha \nabla c)=0$.
Horstmann and Winkler \cite{HorstmannWinkler2005} -- and later many
others -- studied both blow-up and existence of the solutions for
different $\alpha$.

In this paper, we treat the $\nu$-dimensional system with a point
source of strength $Z$. We have the following sufficient condition for
blow-up:
\begin{thm}
  For $\nu\geq 2$, assume that the initial datum $\rho_0$ satisfies
  \begin{eqnarray}
    \label{blowup}
    \gM_\nu^{\nu-2\over\nu}(0)< \dfrac{1}{(\nu-1)2^\nu
    |\bS^{\nu-1}|}\gM_0^{2-\frac{2}{\nu}} - \dfrac{Z}{2(\nu-1)
    |\bS^{\nu-1}|}\gM_0^{1-\frac{2}{\nu}},
  \end{eqnarray}
  where $|\mathbb{S}^{\nu-1}|$ is the volume of $\nu-1$-dimensional
  sphere. Then there is no classical solution that exists for all
  times, i.e., there exists a finite $T>0$ such that
  $\lim_{t\rightarrow T-} \|\rho(\cdot,t)\|_{\infty} =+\infty$.
\end{thm}
This has two immediate consequences:
\begin{description}
\item[No point source]
  Without an external point source, i.e., $Z=0$, the condition
  \eqref{blowup} is a blow-up criterion for the multi-dimensional
  parabolic-elliptic Patlak-Keller-Segel system, i.e., if, initially
  $$
  \gM_\nu^{\nu-2\over\nu}(0)< \dfrac{1}{(\nu-1)2^\nu |\bS^{\nu-1}|}\gM_0^{2-\frac{2}{\nu}},
  $$
  then the solution blows up in finite time. This is actually the
  condition we are looking for, since for $\nu=2$, it is exactly
  Dolbeault's and Perthame's condition \cite{DolbeaultPerthame2004}
  $$
  1<\frac{1}{8\pi}\gM_0.
  $$
\item[With point source and $\nu=2$] In dimension two, the condition
  becomes
  $$
  1+\frac{Z}{4\pi}<\frac{1}{8\pi}\gM_0,
  $$
  which is exactly Wolansky's and Espejo's blow-up condition
  \cite{WolanskyEspejo2013}.
\end{description}

\begin{rem}
It is a interesting question whether the reverse inequality, i.e.,
\begin{eqnarray*}
  \gM_\nu^{\nu-2\over\nu}(0)> \dfrac{1}{(\nu-1)2^\nu
  |\bS^{\nu-1}|}\gM_0^{2-\frac{2}{\nu}} - \dfrac{Z}{2(\nu-1)
  |\bS^{\nu-1}|}\gM_0^{1-\frac{2}{\nu}},
\end{eqnarray*}
would already imply existence of a classical solution. It has an
affirmative answer in dimension two but is open in higher dimensions.
\end{rem}

\section{Proof of the main result}
\begin{proof}
  By using the fundamental solution $\Phi$ of the Poisson equation,
  the system \eqref{equ} can be rewritten into the following form, for
  $\nu\geq 2$,
\begin{equation}
  \label{eq:ks}
  \partial_t\rho=\Delta\rho-\nabla[\rho \mathcal{K}*(\rho-Z\delta)], \qquad x\in\rz^\nu,
\end{equation}
where
$$
\mathcal{K}=\nabla\Phi=-\frac{1}{|\mathbb{S}^{\nu-1}|}\frac{x}{|x|^\nu}.
$$
Multiplication of \eqref{eq:ks} by $|x|^\nu$ and integration gives
\begin{equation}
  \begin{split}
  \label{eq:m3}
  &{\rd\over\rd t}\int_{\rz^\nu} |x|^\nu\rho(x)\rd x\\
  = &2\nu(\nu-1)\int_{\rz^\nu}|x|^{\nu-2} \rho(x)
  \rd x+ {\nu Z\over|\mathbb{S}^{\nu-1}|}\int_{\rz^\nu}\rho(x)\rd x \\
  & -{\nu\over|\mathbb{S}^{\nu-1}|}\int_{\rz^\nu}\int_{\rz^\nu}
  |x|^{\nu-1}{x\over|x|}{x-y \over |x-y|^\nu}\rho(x)\rho(y) \rd x \rd y\\
  =& 2\nu(\nu-1)\int |x|^{\nu-2} \rho(x) \rd x
  + {\nu Z\over|\mathbb{S}^{\nu-1}|}\int_{\rz^\nu}\rho(x)\rd x\\
  &- {\nu\over 2|\mathbb{S}^{\nu-1}|}\int_{\rz^\nu}\int_{\rz^\nu}
  \underbrace{(|x|^{\nu-1}{x\over|x|}-|y|^{\nu-1} {y\over|y|}){x-y
      \over |x-y|^\nu}}_{=:V}\rho(x)\rho(y)\rd x\rd y
\end{split}
\end{equation}
We estimate $V$ from below. To this end we set $r:=|x|$, $s:=|y|$, and
$u:=x\cdot y/(|x||y|)$. Thus
\begin{equation}
  \begin{split}
  \label{v}
  V=&{r^\nu-(r^{\nu-1} s+rs^{\nu-1})u + s^\nu\over (r^2+s^2-2rsu)^{\nu/2}}\\
  =&{(r/s)^{\nu/2}+(s/r)^{\nu/2}-[(r/s)^{\nu/2-1}+(s/r)^{\nu/2-1}]u\over
\{r/s+s/r-2u\}^{\nu/2}}.
\end{split}
\end{equation}
The right hand side of \eqref{v} is monotone increasing in $u$. To see
this we first write $\tau:=r/s$ and remark that the derivative of $V$
with respect $u$ is -- up to an irrelevant non-negative factor --
\begin{equation}
  \begin{split}
   &-\left(\tau^{\frac\nu2-1}+\tau^{-(\frac\nu2-1)}\right)
       (\tau+\frac1\tau-2u)+\nu\left(\tau^{\frac\nu2}+\tau^{-\frac\nu2} -(\tau^{\frac\nu2-1}+\tau^{-(\frac\nu2-1)})u\right)  \\
  \geq & (2-\nu)\left(\tau^{\frac\nu2-1}+\tau^{-(\frac\nu2-1)}\right) +\nu \left(\tau^{\frac\nu2}+\tau^{-\frac\nu2}\right) -\left(\tau+\frac1\tau\right) \left(\tau^{\frac\nu2-1}+\tau^{-(\frac\nu2-1)}\right)\\
  = & (2-\nu)\left(\tau^{\frac\nu2-1}+\tau^{-(\frac\nu2-1)}\right) +(\nu-1) \left(\tau^{\frac\nu2}+\tau^{-\frac\nu2}\right) - \left(\tau^{\frac\nu2-2}+\tau^{-(\frac\nu2-2)}\right)\\
  \geq & (2-\nu)\left(\tau^{\frac\nu2-1}+\tau^{-(\frac\nu2-1)}\right) +(\nu-2) \left(\tau^{\frac\nu2}+\tau^{-\frac\nu2}\right)\geq  0
\end{split}
\end{equation}

where we have used that $\tau^\alpha+\tau^{-\alpha}$ is monotone increasing in $\alpha$ for positive $\alpha$ and $\tau$.
Thus
\begin{eqnarray*}
V&\geq &{(r/s)^{\nu/2}+(s/r)^{\nu/2}+(r/s)^{\nu/2-1}+(s/r)^{\nu/2-1}\over
\{r/s+s/r+2\}^{\nu/2}}\\
&=&{(r/s)^{\tfrac\nu2}+(s/r)^{\tfrac\nu2}+(r/s)^{\tfrac\nu2-1}+(s/r)^{\tfrac\nu2-1}\over
(\sqrt{r/s}+\sqrt{s/r})^\nu}\geq \min_{\tau\in [0,1]} f(\tau)
\end{eqnarray*}
with
$$
f(\tau)={1+\tau^\nu+\tau+\tau^{\nu-1}\over (1+\tau)^\nu}=\frac{1+\tau^{\nu-1}}{(1+\tau)^{\nu-1}}.
$$
Next we prove that $f$ is decreasing and the minimum achieved at $\tau=1$. In fact,
\begin{eqnarray*}
&&f'(\tau)={(\nu-1)\tau^{\nu-2}(1+\tau)^{\nu-1}
-(\nu-1)(1+\tau)^{\nu-2}(1+\tau^{\nu-1})\over(1+\tau)^{2(\nu-1)}}\\
%&=&\frac{\nu-1}{(1+\tau)^{\nu}}[(\nu \tau^{\nu-1}+1+(\nu-1)\tau^{\nu-2})(1+\tau)-\nu(1+\tau^\nu+\tau+\tau^{\nu-1})]\\
&=&\frac{\nu-1}{(1+\tau)^{\nu}}[\tau^{\nu-2}-1]\leq 0.
\end{eqnarray*}
Thus
$$
V\geq f(1) =2^{2-\nu}
$$
and therefore
\begin{equation}
{\rd\over\rd t} \mathfrak{M}_\nu \leq 2\nu(\nu-1)\mathfrak{M}_{\nu-2}-\frac{\nu2^{1-\nu}}{|\mathbb{S}^{\nu-1}|}\gM_0^2+{\nu Z\over|\bS^{\nu-1}|}\gM_0.
\end{equation}
Estimating by H\"older's inequality yields
\begin{equation}
  \label{Hoelder}
  \gM'_\nu\leq 2\nu(\nu-1) \gM^{\nu-2\over\nu}_\nu\gM_0^{2\over\nu}
  -{\nu2^{1-\nu}\over|\bS^{\nu-1}|}\gM^2_0+{\nu Z\over|\bS^{\nu-1}|}\gM_0.
\end{equation}
In particular, we have a shrinking $\nu$-th moment, if the initial moments fulfill
$$
\gM_\nu^{\nu-2\over\nu}(0)< \dfrac{1}{(\nu-1)2^\nu |\bS^{\nu-1}|}\gM_0^{2-\frac{2}{\nu}} - \dfrac{Z}{2(\nu-1) |\bS^{\nu-1}|}\gM_0^{1-\frac{2}{\nu}},
$$
However, a shrinking $\nu$-the moment implies blow-up.
\end{proof}
At this point we would like to remark that the strategy of the proof,
namely multiplication with an appropriate power based on a dimensional
analysis, has been previously used in the context of effective quantum
models and dates back -- at least -- to an unpublished observation of
Benguria to bound the excess charge of atoms. Later, Lieb
\cite{Lieb1984} extended the argument to the quantum case;
Lenzmann and Lewin
\cite{LenzmannLewin2013} used it in the time dependent setting.

In conclusion, we would like to point that our blow-up condition
\eqref{blowup} implies $\eqref{m2m0}$. To see this, we first note that
$\gM_2\leq \gM_0^{\nu-2\over\nu} \gM_\nu^{2\over\nu}$ by
interpolation. Estimating the right hand side by using
$\eqref{blowup}$ gives $\eqref{m2m0}$ with the extra bonus of a
definite constant instead of an uncontrolled one. Finally, note that
this is only necessary when $\nu>2$, since Inequality $\eqref{m2m0}$
is an empty statement in dimension two whereas Inequality
$\eqref{blowup}$ remains meaningful.

{\small \textit{Acknowledgment:} We thank Georgios Psaradakis for
  critical reading of the manuscript. We acknowledge support by the
  Deutsche Forschungsgemeinschaft through the grants CH 955/4-1 and SI
  348/15-1.}

\def\cprime{$'$}

\end{document}